\documentclass{article}
\usepackage{amssymb}
\usepackage{amsmath}

\let\<\langle
\let\>\rangle
\font\bigcmsy=cmsy10.pk scaled 2000
\def\bigtimes{\mathop{\,\vrule width0pt depth2pt height8pt
            \smash{\lower2pt\hbox{\bigcmsy\char'002}}\,}\limits}

\begin{document}

\begin{center}
\Large{\textbf{Generalization of order separability for free groups.}}
\end{center}

\begin{center}
\textbf{Vladimir V. Yedynak}
\end{center}

\begin{abstract}
In this work the author studies the property close to property of order separability.


\end{abstract}

\section{Introduction.}

Definition 1. We say that a group $G$ is \textsl{order separable}
if, for each elements $g, h$ of $G$ such that $g$ and $h^{\pm1}$ are
not conjugate, there exists a homomorphism of $G$ onto a finite
group such that the orders of the images of $g$ and $h$ are
different.

In [1] it is proved that free groups are order separable. This result was generalized in [2] where it was shown that each free group $F$ is omnipotent that is for an arbitrary finite set of elements $g_1,..., g_n$ of $F$ such that for each $i\neq j$ elements $g_i$ and $g_j$ have no nontrivial conjugate powers there exists a constant $k$ such that for each ordered sequence of positive integers $l_1,..., l_n$ there exists a homomorphism $\varphi$ of $F$ onto a finite group such that the order of $\varphi(g_i)$ is $kl_i$.

It is known also that the property of order separability is inherited by free products [3]. The following theorem generalizes the property of order separability for free groups.

\textbf{Theorem.} \textsl{Let $u_1, u_2$ be the elements of a free group $F$, which does not belong to conjugate cyclic subgroups. Then for each prime number $p$ and for each integer $m$ there exists a homomorphism $\varphi$ of $F$ onto a finite $p$-group such that $u_1$ and $u_2$ do not belong to the kernel of $\varphi$ and $\mid\varphi(u_1)\mid
/\mid\varphi(u_2)\mid=p^n.$}

\section{Notations and definitions.}

We shall use the correspondence between the actions of a free group
$F(x, y)$ with basis $x, y$ and the graphs satisfying the following
properties:

1) for each vertex $p$ of a graph there exist exactly one edge with
label $x$ that goes away from this vertex and exactly one edge
with label $y$ that goes into $p$;

2) we consider that for each labelled edge there exists the corresponding inverse edge and each two edges with labels are not mutually inverse;

Definition 2. We say that the graph is the \textsl{action graph} of
a free group $F=F(x, y)$ if it satisfies the properties 1), 2).
We shall consider additionally that all labelled edges in this graph
are positively oriented and set the orientation of this graph.

If $\varphi$ is a homomorphism of a group $F$ then the Cayley graph
of the group $\varphi(F)$ with the generating set $\{\varphi(x),
\varphi(y)\}$ is the action graph for $F$ (we identify labels
$\varphi(x), \varphi(y)$ with $x, y$ correspondingly in this graph).

Let $\Gamma$ be the action graph of the free group $F$. Fix a
natural number $n$, vertex $p$ of the graph $\Gamma$ and element
$z\in\{x^{\pm1}, y^{\pm1}\}$ of the group $F$. Consider $n$ copies
of the graph $\Gamma: \Delta_1,..., \Delta_n$. Let $p_i$ be the
vertex of the graph $\Delta_i$ corresponding to the vertex $p$ of
the graph $\Gamma$. If $z\in\{x, y\}$ then denote by $q_i$ the
vertex such that there exists the edge $f_i$ with label $z$ which
goes into $q_i$ from $p_i$. If $z\in\{x^{-1}, y^{-1}\}$ then denote
by $q_i$ the vertex such that there exists the edge $f_i$ with label
$z^{-1}$ which goes goes away from $q_i$ into $p_i$. We construct
the graph $\Delta=\gamma_n(\Gamma; p; z)$ from graphs $\Delta_i,
1\leqslant i\leqslant n$ deleting edges $f_i$ and connecting
vertexes $p_i$ and $q_{i+1}$ by an edge whose label equals either
$z$ when $z\in\{x, y\}$ or $z^{-1}$ when $z\in\{x^{-1}, y^{-1}\}$
(indices are modulo $n$). If $z\in\{x, y\}$ then this new edge goes
away from the vertex $p_i$. If $z\in\{x^{-1}, y^{-1}\}$ then this
new edge goes into the vertex $p_i$. The graph $\Delta$ is the
action graph of the group $F$.

If $S$ is a path in a graph then $\alpha(S), \omega(S)$ are the
beginning and the end of $S$ correspondingly. If $S=e_1...e_n$ is a
path in the action graph of the group $F$, then Lab$(S)= \
$Lab$(e_1)'...$Lab$(e_k)'$ is the label of the path $S$, where
Lab$(e_i)'= \ $Lab$(e_i)$ is a label of edge $e_i$ in case $e_i$ is
positively oriented and if the edge $e_i$ is negatively oriented
then Lab$(e_i)'=$Lab$(e_i')^{-1}$, where edges $e_i$ è $e_i'$ are
mutually inverse.

Definition 3. Suppose we have two action graphs of the group $F$
--- $\Gamma$ and $\gamma_n(\Gamma; q; z)$. If $p$ is a vertex
of $\Gamma$ then $p^i$ is a vertex in the graph $\gamma_n(\Gamma; q;
z)$ which belongs to $i$-th copy of the graph $\Gamma$ and
corresponds to the vertex $p$. If $S$ is a path in $\Gamma$ then
$S^i$ is the path in $\gamma_n(\Gamma; q; z)$ which goes from the
vertex $\alpha(S)^i$ and whose label is equal to the label of the
path $S$.

Definition 4. Let $u$ be a cyclically reduced element of the free
group $F$. If $S=e_1...e_n$ is a closed path without returnings
whose label equals $u^k$ in the action graph of $F$ then we say that
$S$ is the $u$-cycle whose length equals $k$ (we consider that there
exists exactly one subpath of the path $S$ ending in $\omega(S)$
whose label equals $u$).

\section{Proof of theorem.}

Fix elements $u_1, u_2$ of $F$. We consider that these elements are cyclically reduced.

Suppose that $u_i=u_i'^{p^{m_i}t_i}, i=1, 2,$ where $p\nmid t_i$.
Denote by $v_i$ the elements $u_i^{t_i}, i=1, 2$. Without loss of
generality we may assume that $m_1 \geqslant m_2$. If there exists a
homomorphism $\psi$ of $F$ onto a finite group such that
$\mid\psi(v_1)\mid/\mid\psi(v_2)\mid=p^{n+m_1-m_2}$ then
$\mid\psi(u_1)\mid/\mid\psi(u_2)\mid=p^n$. Hence we may assume that
the subgroups $\<u_1\>$ and $\<u_2\>$ are $p'$-isolated.

There exists a homomorphism $\varphi$ of $F$ onto a finite $p$-group
such that the images of the elements $u_1, u_2$ are nonunit
$p$-elements [4]. Since all $p'$-isolated cyclic subgroups of free
groups are separable in the class of finite $p$-groups we may also
assume that all $u_1$- and $u_2$-cycles in the Cayley graph $\Gamma$
of the group $\varphi(G)$ with the set of generators $x=\varphi(x'),
y=\varphi(y')$ are simple. Without loss of generality we may
consider that $|\varphi(u_1)|=p^k, |\varphi(u_2)|=p^{k+l},
k\geqslant1$. Put $\Gamma_{-1}=\Gamma$. Let
$u_1=y_0^{\varepsilon_1}...y_{k-1}^{\varepsilon_1}, y_i\in\{x, y\},
\varepsilon_i\in\{-1, 1\}$ be the reduced form of the element $u_1$
in the basis $x, y$. Fix an arbitrary vertex $q$ in the graph
$\Gamma_{-1}$. For $i>-1$ we shall define by the induction the graph
$\Gamma_i=\gamma_p(\Gamma_{i-1}; q_{i-1}; y_{i'}^{\varepsilon_{i'}})$
and the path $S_i$ in $\Gamma_i$ whose length equals $i+1$, where
$i'$ is the remainder from a division of $i$ on $k$. If $i>-1$ then $q_i$
is the vertex of the graph $\Gamma_i$ which is the end point of the
path $S_i$. The vertex $q_{-1}$ is equal to the vertex $q$. If
$i>-1$ we define a path $S_i$ in the graph $\Gamma_i$ in the
following way. If $i=0$ then $S_0$ is the first edge of the
$u_1$-cycle which goes from the vertex $q_{-1}^1$ or its inverse. In case $i>0$ we
define the path $S_i$ as $S_{i-1}^1\cup f_i$ where $f_i$ is the edge
of the graph $\Gamma_i$ one of whose endpoints coincides with
$\omega(S_{i-1}^1)$ and the label of $f_i$ or its inverse equals $y_{i'}$. Also if
$\varepsilon_i=1$ then $f_i$ is positively oriented and has a label. If $\varepsilon_i=-1$ then $f_i$ is negatively oriented. It is easy to notice that the length of
each $u_1$- or $u_2$-cycle in $\Gamma_i$ is the power of $p$. Since
all $u_1$- and $u_2$-cycles are simple in the graph $\Gamma_{-1}$
this condition is also held for $\Gamma_i$ for each $i$. Besides for
each $i$ the graph $\Gamma_i$ contains the maximal $u_1$-cycle which
contains the path $S_i$.

Suppose that there exists $i$ such that the following conditions are
true. In the graph $\Gamma_j$ each maximal $u_2$-cycle contains
the path $S_j$ for all $j\leqslant i$ but not for $j=i+1$ (if
$i+1=0$ we simply consider that in the graph $\Gamma_0$ not all
maximal $u_2$-cycles contains $S_0$). Notice that if
$j\leqslant i$ then in the graph $\Gamma_j$ the length of the
maximal $u_i$-cycle coincide with $|\varphi(u_i)|p^{j+1}, i=1, 2$.
In the graph $\Gamma_{i+1}$ the length of the maximal $u_1$-cycle
equals $\mid\varphi(u_1)\mid p^{i+1}$. Lets find the length of the
maximal $u_2$-cycle in $\Gamma_{i+1}$. Fix a vertex $r$ of
$\Gamma_i$. Then the length of the $u_2$-cycle which goes away from
$r$ is equal to $|\varphi(u_2)|p^k, k\leqslant i$. Then for each
$l=1,..., p$ in the graph $\Gamma_{i+1}$ the length of the
$u_2$-cycle that goes away from the vertex $r_l$ is not more than
$|\varphi(u_2)|p^{k+1}$. Suppose that the $u_2$-cycle $T$ of
$\Gamma_i$ starting from $r$ is maximal, that is its length coincide
with $|\varphi(u_2)|p^i$. The path $S_i$ is contained in $T$.
The $u_2$-cycle $T^1$ of $\Gamma_{i+1}$ contains the path $S_i^1$.
From the condition on $i$ and the simplicity of $u_2$-cycles of
$\Gamma_{i+1}$ we may deduce that $T^1$ does not contain the edge
$f_{i+1}$. This $u_2$-cycle does not also contain the edge
connecting the first and the last copies of the graph $\Gamma_i$
because otherwise $T$ would have the self-intersection in the vertex
$\omega(S_{i-1})^1$. Hence the length of $T^1$ coincide with the
length of $T$. It means that the length of the maximal $u_2$-cycle
in $\Gamma_{i+1}$ coincide with $\mid\varphi(u_2)\mid p^i$.
Since $u_1$ and $u_2$ does not belong to the conjugate cyclic
subgroups then there exists $i$ that satisfies the conditions
mentioned above.

The number of vertices in the graph $\Gamma_{i+1}$ equals
$|\varphi(F)|p^{i+2}$. So there exists a homomorphism $\varphi_1$ of
$F$ onto a finite $p$-group such that in the Cayley graph of the
group $\varphi_1(F)$ with the set of generators $\varphi_1(x'),
\varphi_1(y')$ all $u_1$- and $u_2$-cycles are simple. Besides
$|\varphi_1(u_2)|/|\varphi_1(u_1)|=|\varphi(u_2)|/|\varphi(u_1)|*1/p$.

If we construct graphs $\Gamma_i$ using the element $u_2$ instead of
$u_1$ we obtain the homomorphism $\varphi_2$ that satisfies the same
conditions concerning to the $u_1$- and $u_2$-cycles in the graph
$\varphi_2(F)$ and the homomorphism $\varphi_1$ and
$|\varphi_2(u_2)|/|\varphi_2(u_1)|=|\varphi(u_2)|/|\varphi(u_1)|*p$.
In order to obtain the homomorphism of $F$ onto a finite $p$-group
such that the ratio of orders of the images of $u_1$ and $u_2$
sufficient $p^n$ it is enough to apply these procedures several times.
The theorem is proved.

\begin{center}
\large{Acknowledgements.}
\end{center}

The author thanks A. A. Klyachko for valuable comments and discussions.

\begin{center}
\large{References.}
\end{center}

1. \textsl{Klyachko A. A.} Equations over groups, quasivarieties,
and a residual property of a free group // J. Group Theory. 1999.
\textbf{2}. 319--327.

2. \textsl{Wise, Daniel T.} Subgroup separability of graphs of free groups with cyclic edge groups.
Q. J. Math. 51, No.1, 107-129 (2000). [ISSN 0033-5606; ISSN 1464-3847]

3. Yedynak V. V. Order separability. \textsl{Bulletin of MSU}
\textbf{3}: 53-56.

4. Kargapolov, M. I., Merzlyakov, Yu. I. (1977). Foundations of group
theory. \textsl{Nauka}.

\end{document}